\documentclass[10pt]{siamltex}
\usepackage{cite,epsfig,graphicx}

\def\eqnok#1{\/{\rm (\ref{eq:#1})}}
\def\begeq{\begin{equation}}
\def\endeq{\end{equation}}
\def\ie{{\em i.\ e.\ \/}}

\def\calp{{{\cal P}}}
\def\calk{{{\cal K}}}

\title{
Newton-Krylov solvers for time-steppers
\thanks{Version of \today.}
}

\author{C. T. Kelley\footnotemark[2]
\and I. G. Kevrekidis\footnotemark[3]\ \footnotemark[4]
\and L. Qiao\footnotemark[3]}

\begin{document}
\maketitle
\renewcommand{\thefootnote}{\fnsymbol{footnote}}
\footnotetext[2]{North Carolina State University, Center for
Research in Scientific Computation and Department of Mathematics,
Box 8205, Raleigh, N.C. 27695-8205, USA ({Tim\_Kelley@ncsu.edu}).
This research was supported in part by National Science Foundation
grants DMS-0070641 and DMS-0209695.} \footnotetext[3]{Department
of Chemical Engineering, Princeton University, Princeton, NJ 08544
({yannis@princeton.edu}). Work supported in part by AFOSR and an
NSF/ITR grant.} \footnotetext[4]{PACM and Mathematics, Princeton
University.}\renewcommand{\thefootnote}{\arabic{footnote}}

\begin{abstract}

We study how the Newton-GMRES iteration
can enable dynamic simulators (time-steppers) to
perform fixed-point and path-following computations.
For a class of dissipative problems, whose dynamics are characterized
by a slow manifold, the Jacobian matrices in
such computations are compact perturbations of the identity.
We examine the number of GMRES iterations required
for each nonlinear iteration as a function of the dimension of the slow subspace
and the time-stepper reporting horizon.
In a path-following computation, only a small number (one or two) of
additional GMRES iterations is required.
%
%
\end{abstract}



\section{Introduction}
\label{sec:intro}

In studying the dynamic behavior of evolution equations, \begeq
\label{eq:fode}
du/dt=f(u;\lambda),
\endeq
a computational modeler typically chooses between two paths:
the first is developing a dynamic simulator for the process; the
second is developing algorithms to locate particular features of
the long-term dynamics of the process, such as steady states or
limit cycles.
The first path typically gives rise to initial value problems, and the
corresponding codes are dynamic simulators - we will call them ``time-steppers"
since, given the state of the system at a moment in time they produce (an
approximation of) the state of the system at a later moment.
The second path typically gives rise to fixed point algorithms for solving
coupled nonlinear algebraic equations; these may be the steady state equations
themselves, or an augmented set arising in a continuation/bifurcation context.
The Recursive Projection Method of Shroff and Keller
\cite{herbrpm} (see also the Adaptive Condensation of Jarausch and
Mackens \cite{J&M1987a,J&M1987b} for symmetric problems and the Newton-Picard work of Lust {\it et al.}\cite{lust}
 ) is an example of an algorithm that, in some sense, connects the two paths.
The main idea is to construct a computational superstructure that
designs and combines several calls to an existing (``legacy")
time-stepper, effectively turning it into a fixed point solver for
\[  u - \Phi_T(u;\lambda) = 0, \]
where $\Phi_T$ is the result of integration of  \eqnok{fode} with
initial condition $u$ for time T (see the review of Tuckerman and Barkley \cite{tuckerman}).
Here the action of the linearization of the time-stepper is estimated
in a matrix-free fashion by the integration of appropriately chosen nearby initial
conditions.

Over the last few years,
this matrix-free computational enabling technology has found
new applications in the field of multiscale computations.
In current modelling practice, dynamic models are often constructed
at a microscopic/stochastic level of description (e.g. molecular
dynamics, kinetic Monte Carlo or Lattice Boltzmann codes); the
closure required to obtain explicit, macroscopic ``effective"
equations is not available in closed form.
Algorithms like RPM can then be ``wrapped around" appropriately
initialized ensembles of short bursts of microscopic simulations,
and solve the effective equation without ever obtaining it in
closed form.
Several examples of such ``coarse equation-free" computation have
now been explored, and the mathematical underpinnings of the
approach are being extensively studied
\cite{coarse1,coarse2,coarse3,coarse4}.

In this paper we study the Newton-GMRES iteration as a computational
``wrapper" around a legacy time-stepper.
This wrapper enables the computation and continuation of fixed points of the time-$T$
map of the time-stepper ({\it i.e.} steady states of the corresponding dynamical equations).
It is useful, for purposes of discussion, to consider that the time-stepper is
available as an input-output black box  (an executable) which cannot be modified.

The paper is organized as follows:
We first review certain properties of GMRES in the context of problems whose
linearization
is a compact perturbation of the identity.
We argue that time-steppers for a class of dissipative problems, whose long-term
dynamics lie on a low-dimensional, slow manifold, may fit this description; we then
proceed to examine Newton-GMRES convergence and number of iterations
for fixed point computation of such time-steppers in a continuation context.
Our numerical examples are the Chandrasekhar H-equation
\cite{chand} and the time-stepper of a discretization of a reaction-diffusion problem
known to possess a low-dimensional inertial manifold
\cite{chaffee,temam}.
We illustrate the effect, on the GMRES iterations, of the time-stepper reporting horizon
both for fixed point and for pseudo-arclength continuation computations.
We conclude with a brief discussion and thoughts about extensions of this approach.

\section{Convergence Analysis}
\label{sec:converge}

\subsection{GMRES Preliminaries}
\label{subsec:prelim}

We will use Newton-GMRES to solve the nonlinear fixed point
equation \begeq \label{eq:ftstepper} u - \Phi_T(u;\lambda) \equiv
F(u)=0.
\endeq
In what follows we will work in $R^N$ and use the usual Euclidean
norm; the idea is that equation \eqnok{ftstepper} arises from
integrating a set of ODEs, possibly from the discretization of a
partial differential equation on a given mesh.
We will also assume (and the consequences of this will become
apparent below) that the long-term dynamics of the discretized PDE
occur on an attracting ``slow manifold" of low dimension $p < < N$; $p$ will be
assumed independent of mesh refinement ({\it i.e.} it will remain
constant as N increases.)
%


GMRES is an iterative method for solving linear systems $A x = b$ in $R^N$.
The $k$th GMRES iteration minimizes the residual $r$ over
$x_0 + \calk_k$, where $x_0$ is the initial iterate and
$\calk_k$ is the the $k$th Krylov subspace
\[
\calk_k = \mbox{span} \{ r_0, A r_0, \dots , A^{k-1} r_0 \}.
\]
A consequence \cite{gmres,ctk:roots}
of the minimization is that
\begeq
\label{eq:gmresmin}
\| r_k \| = \min_{p \in \calp_k } \| p(A) r_0 \|
\endeq
where $\calp_k$ is the space of $k$th degree residual polynomials, {\it i.e.}
polynomials of degree $k$ such that $p(0) = 1$.

We will apply \eqnok{gmresmin} to a special class of problems where
\begeq
\label{eq:class}
A = I - K + E.
\endeq
In \eqnok{class}
\begin{itemize}
\item $I - K$ is nonsingular,
\item $K = P_D K P_D$, where $P_D$ is an orthogonal projection onto a space $D$ of
dimension $p << N$, and
\item $E$ is a matrix with small norm.
\end{itemize}

We will analyze the performance of GMRES in a way different from
the eigenvalue-based approach used for diagonalizable matrices
\cite{kerksaad,ctk:roots,anne,trefbau}. For the class of problems
of interest here, we can prove a convergence result directly using
methods similar to those in
\cite{ctk:jordan,ctk:xue4,ctk:xue3,ctk:ferng1}. The result in this
paper has stronger and sharper convergence rates.
This result carries through in the infinite-dimensional
case also, using the $L^2$ norm for the corresponding function spaces.


\begin{theorem}
\label{th:conv}
Let $A$ be given by \eqnok{class}. Then there exists $C$ such
that for all $m \ge 1$,
\begeq
\label{eq:residualk}
\| r_{m (p+1)} \| \le C^m \| E \|^m.
\endeq
\end{theorem}

\begin{proof}
Let $p_C$ be the characteristic polynomial of $I - K$. Since
$D$ has dimension $p$, $p_C$ has degree $p+1$. Clearly
\[
p_C(A) = p_C(I - K) + \Delta = \Delta.
\]
by the Cayley-Hamilton theorem. Moreover, there is $C > 0$ such that
\[
\| \Delta \| \le C \| E \|.
\]
In fact, if
\[
p(z) = 1 + \sum_{k=1}^{p+1} \gamma_k z^k
\]
then we may use
\[
C = \sum_{k=1}^{p+1} k | \gamma_k | \| I - P_D K P_D \|^{k-1}
+ O(\| E \|^2).
\]
Hence
\[
\| p(A) \| \le C \| E \|.
\]
This is \eqnok{residualk} for $m = 1$.

Define
\[
p(z) = p_C(z)/p_C(0).
\]
Clearly $p \in \calp_{p+1}$. Hence, by \eqnok{gmresmin}
\begeq
\label{eq:terminate}
\| r_{m(p+1)} \| \le \| p^m(A) r_0 \| \le C^m \| E \|^m \| r_0 \|,
\endeq
as asserted.
\end{proof}

The estimate in \eqnok{terminate} does not depend on the eigenvalues
of $A$, nor is there any requirement that $A$ be normal or even
diagonalizable.

As is standard, a GMRES iteration is terminated when
\[
\| r_k \| \le \eta \| r_0 \|,
\]
where $\eta$ is an user-defined parameter. In the context of this
paper $\| E \|$ is well below the termination tolerance $\tau$, so
we conclude that the iteration will terminate in at most $p+1$ iterations.
In the general case, of course, more cycles could be required.

\subsection{Time-Steppers and Steady State Solutions}
\label{subsec:steady}

Let $\Phi(u)$ denote the output of the time-stepper with time step
$T$ and initial data $u$; we have dropped the subscript $T$ of
$\Phi$ for convenience.
We seek to solve
\begeq
\label{eq:ftime}
F(u) \equiv u - \Phi(u) = 0,
\endeq
to find a steady state solution.
Consider the structure of the eigenvalues
$\mu_i = \exp(\sigma_i T)$ of $\Phi_u$ given the
structure of the eigenvalues $\sigma_i$
of the linearization of the original problem.
We assume that there exists a significant gap between a few
$\sigma_i$ close to zero ($p$ of them, to be exact) and the
remaining large negative $\sigma_j$.
The $p$ small norm $\sigma_i$ do not have to be stable - they can
also be slightly unstable, so that the corresponding $\mu_i$ are
close to the unit circle.

We then assume that $T$ is large enough
so that there is a space $D$ of low dimension $p$, so that
\[
\Phi(u) = P_D \Phi(P_D u) + E(u)
\]
and $E$ and its Jacobian $E_u$ are small.
For a problem where some of the $p$ small $\sigma_i$ are unstable,
this large time $T$ should not be so large that the unstable modes
will cause the trajectory to move significantly away from the
fixed point. In such cases the proper choice of $T$ is a somewhat delicate
matter \cite{herbrpm,coarse2}.

The equation for the Newton step from a current iterate $u_c$ is
\begeq \label{eq:nstep} F_u(u_c) s = s - \Phi_u(u_c) s = -F(u_c),
\endeq
here $F_u$ and $\Phi_u$ are the Jacobians of $F$ and $\Phi$.

A Newton-GMRES \cite{ctk:roots,brown/saad90,brown/saad94} method
solves \eqnok{nstep} with GMRES, terminating the linear
(or inner) iteration when
\[
\| F_u(u_c) s \| \le \eta_c \| F(u_c) \|
\]
where $\eta_c$ may be changed as the outer (or nonlinear)
iteration progresses
\cite{ctk:roots,homerstan,ctk:newton}. For the application here,
we assume that $\| E_u(u) \|$ is much less than any choice of
$\eta$ we make during the iteration.

The linear system \eqnok{nstep} fits exactly in to the paradigm of
\S~\ref{subsec:prelim} with $K=P_D \Phi_u(P_D u) P_D$ and $E =
E_u(u)$. Hence we conclude that a GMRES iteration for
\eqnok{nstep} will take at most $p+1$ iterations to drive the
residual to $O(\| E \|)$.

\subsection{Time-Steppers and Continuation}
\label{subsec:continuation}

In a continuation context $\Phi$ depends on a parameter $\lambda$. As is
standard \cite{herb} we add an additional arclength parameter $s$ to obtain
the augmented system \begeq \label{eq:augment} G(u,\lambda,s) =
\left( \begin{array}{c}
F(u(s),\lambda(s)) \\
{\dot u}^T (u - u_0)+ {\dot \lambda} (\lambda - \lambda_0) - (s -
s_0)
\end{array}
\right).
\endeq
In \eqnok{augment}, $\dot u$ and $\dot \lambda$ are approximations
to $du/ds$ and $d \lambda/ds$, which can be obtained in several
ways \cite{herb}. In the calculations reported in
\S~\ref{sec:numbers} we used the estimate of the slope given by
the last two points, $(u_0,\lambda_0)$ and $(u_{-1},\lambda_{-1})$
computed on the branch.


$G$ is defined on $R^{N+1}$ and
\[
G_{u,\lambda} (u,\lambda) =
\left( \begin{array}{cc}
F_u & F_\lambda \\
{\dot u^T} & {\dot \lambda}
\end{array}
\right),
\]
We seek to show that $G_{u,\lambda}$ also fits into our paradigm, with $p$
replaced by at most $p+2$.
Hence, at most two additional linear iterations will be needed for
each Newton iteration of the augmented system.

We use \eqnok{ftime} to obtain
\[
G_{u,\lambda} = \left( \begin{array}{cc}
I - P_D \Phi_u(P_D u;\lambda)P_D  & - P_D \Phi_\lambda \\
{\dot u^T}  & {\dot \lambda}
\end{array}
\right) +
\left( \begin{array}{cc}
-E_u & -E_\lambda \\
0 & 0
\end{array}
\right).
\]

Let
\[
{\cal A} = \left( \begin{array}{cc}
I & 0 \\
0 & 1
\end{array} \right ) - {\cal K} + {\cal E}
\]
where,
\[
{\cal K} = \left( \begin{array}{cc}
P_D \Phi_u P_D & P_D \Phi_\lambda \\
- {\dot u}^T  & 1 - {\dot \lambda}
\end{array}
\right),
\]
and
\[
{\cal E} =
\left( \begin{array}{cc}
-E_u & -E_\lambda \\
0 & 0
\end{array}
\right).
\]

This fits the paradigm of \S~\ref{subsec:prelim}. To see that, let

\[
{\cal D} = \mbox{span} \left( \left( \begin{array}{c} D \\ 0
\end{array} \right), \left( \begin{array}{c} {\dot u}  \\ 0
\end{array} \right), \left( \begin{array}{c} 0 \\  1 \end{array}
\right) \right) \subset R^{N+1},
\]

Clearly the range of ${\cal K}$
\[
R({\cal K}) = \mbox{span} \left( \left( \begin{array}{c} D \\ 0
\end{array} \right), \left( \begin{array}{c} 0 \\  1 \end{array}
\right) \right) \subset {\cal D}.
\]
To apply the results from the previous section, however, we need
$\calk = P_{\cal D} \calk P_{\cal D}$, where $P_{\cal D}$
is the orthogonal projector onto $\cal D$.  This is why the extra dimension
$(\dot u, 0)^T$ is required.

To see this, let
$y=(u,\mu)\subset R^{N+1}$ be orthogonal to $\cal D$. This means that
$y=(\omega,0)$, where $\omega$ is orthogonal to both $D$ and
$\dot u$. Clearly
\[
{\cal K}y = \left( \begin{array}{c}
P_D \Phi_u P_D \omega \\
-{\dot u}^T \omega
\end{array}
\right) =\left( \begin{array}{c} 0\\0\end{array} \right),\] so
${\cal K} ( I-P_{\cal D})=0$. Since
$R({\cal K}) \subset {\cal D}$, we have
$P_{\cal D}{\cal K}={\cal K}$. Summarizing
\[
\calk = P_{\cal D} \calk P_{\cal D}
\]
and
\[
\mbox{dim}({\cal D}) \le p+2.
\]

The dimension of ${\cal D}$ can be taken $p+1$ if ${\dot u}$ is nearly in
the range of $P_D$, \ie
\begeq
\label{eq:dotugoal}
\| (I - P_D) {\dot u} \| = O( \| E_u \| + \| E_\lambda \|  ).
\endeq
In this case, we can let
\[
{\cal D} = \mbox{span} \left( \left( \begin{array}{c} D \\ 0
\end{array} \right), \left( \begin{array}{c} 0 \\  1 \end{array}
\right) \right) \subset R^{N+1},
\]
\[
{\cal K} = \left( \begin{array}{cc}
P_D \Phi_u P_D & P_D \Phi_\lambda \\
- (P_D {\dot u})^T  & 1 - {\dot \lambda}
\end{array}
\right),
\]
and
\[
{\cal E} =
\left( \begin{array}{cc}
-E_u & -E_\lambda \\
((I - P_D) {\dot u} )^T & 0
\end{array}
\right).
\]

If $F_u$ is well conditioned and $| {\dot \lambda} |$
and $\| F_\lambda \|$ are
not too large, then \eqnok{dotugoal} holds.
To see this we differentiate $F = 0$ and obtain
\[
F_u {\dot u} + F_\lambda {\dot \lambda} = 0,
\]
which implies, if $F_u$ is nonsingular, that
\begeq
\label{eq:dotuform}
{\dot u} = - F_u^{-1} F_\lambda {\dot \lambda}.
\endeq

Our assumptions are that $\| E_u \|$ and $\| E_\lambda \|$ are
much smaller than $\| F_u \|$ and $\| F_\lambda \|$. Hence, if $F_u$
is well conditioned, the Banach lemma implies that
\[
F_u^{-1} = (I - P_D \Phi_u P_D)^{-1} + O(\| E_u \|).
\]
Moreover,
\[
F_\lambda =-P_D \Phi_\lambda + O(\| E_\lambda \|) .
\]
We incorporate this into \eqnok{dotuform} to obtain \begeq
\label{eq:dotu2} {\dot u} =  (I - P_D \Phi_u P_D)^{-1} P_D
\Phi_\lambda {\dot \lambda} + O\left( | {\dot \lambda} | ( \| E_u
\|  \| F_\lambda \| + \| E_\lambda \| ) \right),
\endeq
which implies \eqnok{dotugoal} if $| \dot \lambda |$ and
$\| F_\lambda \|$ are $O(1)$.

Summarizing, if $\| E_u \|$ and $\|E_\lambda\|$ are much smaller
that $\| F_u \|$ and $\|F_\lambda\|$ respectively,
$| \dot \lambda |$ and $\| F_\lambda \|$ are $O(1)$,
and if $F_u$ is well conditioned, then \eqnok{dotugoal} holds.
Near folds and bifurcations, $F_u$ is singular, and in those
cases we may need to require that the dimension of $\cal D$
be $p+2$.

\section{Numerical Results}
\label{sec:numbers}

All the computations were done using MATLAB version 6.0 Release 12
on a PC with a Pentium4 2.53GHz CPU.

\subsection{The H-equation}

As a first example, we solved a problem which does not arise
in a time-stepper context, but for which the Jacobian fits
our paradigm \cite{twm68}. The solution and continuation problem for a 100-node
midpoint rule discretization of the Chandrasekhar H-equation
\[
F(x)_i=x_i-\left(1-\frac{c}{2N}\sum_{j=1}^{N}\frac{\mu_ix_j}{\mu_i +\mu_j}
\right)^{-1}
\]

 were obtained with a Newton-GMRES solver
{\tt nsoli} from \cite{ctk:newton}. We set the
relative and absolute tolerances in the solver to
$10^{-12}$.

\begin{figure}[htbp]
\begin{center}
\includegraphics[width=4.5in]{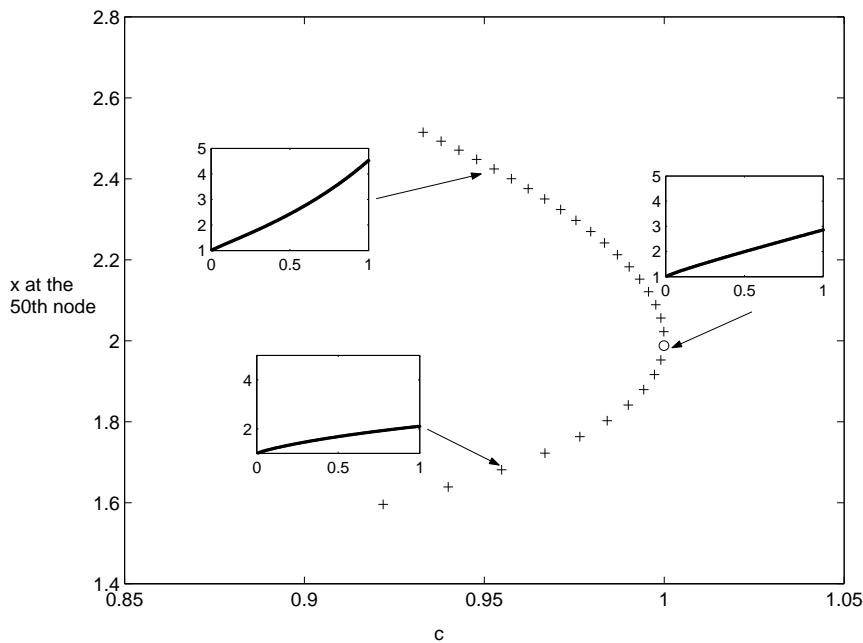}
\caption{Bifurcation diagram of the H-equation for $c$ close to 1. The
parameter $c$ is equal to 0.9999179 at the point marked by a circle.
} \label{hbdiag}
\end{center}
\end{figure}

Figure.~\ref{hbdiag} is a bifurcation diagram with respect to the
parameter c, showing a turning point at $c =1$.
For all values of $c$ shown in this continuation, the eigenvalues
of the linearization of the solution are close to 1; only a single
one of them changes in marked way, ranging from -0.9 on the upper
branch through zero at the turning point to 0.5 on the lower
branch.

\begin{figure}[htbp]
\begin{center}
\includegraphics[width=4.5in]{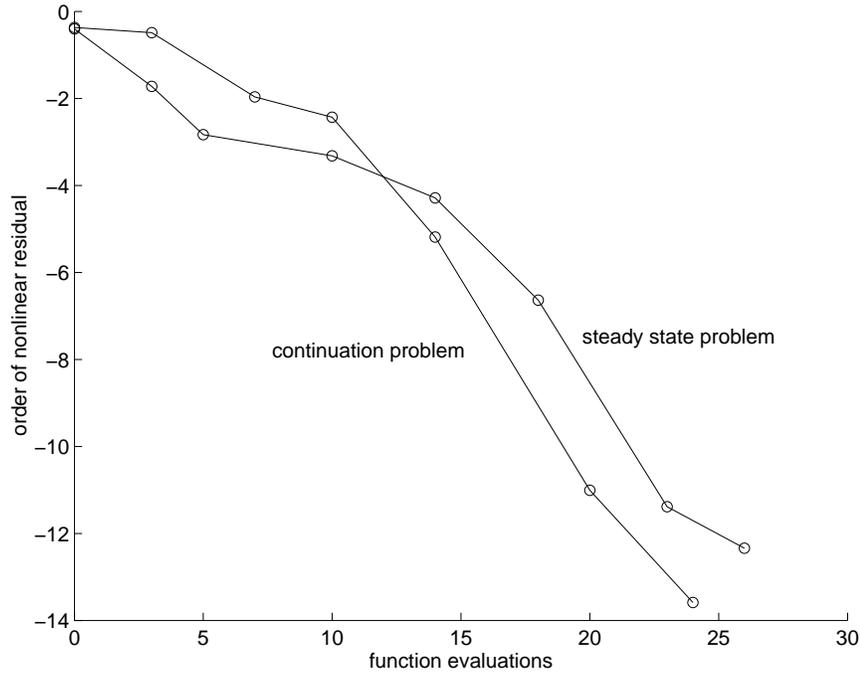}
\caption{Convergence plot for the steady state and continuation
problem of H-equation} \label{fig4}
\end{center}
\end{figure}

Figure.~\ref{fig4} shows the GMRES performance for a
representative value of $c$ ($0.99991$ \\ $79$), close to the
turning point, marked by a circle in Figure~\ref{hbdiag}. We
report on computations from the continuation itself, where the
initial iterate was the standard linear predictor, and from a
second stand-alone computation where the initial iterate was the
solution plus a perturbation function $0.05 \sin(x)$.
In Figure~\ref{fig4} we plot the convergence history for Newton-GMRES
in the two cases. The convergence rates and
the costs of the solves are roughly the same.

The first column of Table.~\ref{heig} shows the ten eigenvalues
farthest from $1$ for the linearization of each of the two problems.
In the continuation case, there is one more eigenvalue far away
from $1$ than in the stand-alone case. In both cases the linearization
clearly fits the pattern of a compact perturbation of the identity.

\begin{table}[htbp]
\begin{center}
\begin{tabular}{|c|c|} \hline
Eigenvalues of steady state problem & Eigenvalues of continuation
problem\\\hline

   0.0207265 & -5.2022448\\
   0.9424114 & 5.2008250\\
   0.9900391 & 0.9746094\\
   0.9974043 & 0.9828252\\
   0.9992541 & 0.9977360\\
   0.9998164 & 0.9991859\\
   0.9999612 & 0.9998098\\
   0.9999926 & 0.9999577\\
   0.9999987 & 0.9999919\\
   0.9999998 & 0.9999985\\\hline

\end{tabular}
\end{center}
\caption{Eigenvalues for the linearized steady state and
continuation problems of the H-equation at $c=0.9999179$}
\label{heig}
\end{table}

\subsection{The Chafee-Infante reaction diffusion problem}

We then solved the steady state and continuation problems for a
discretization of a dissipative reaction-diffusion PDE in one
dimension, the so-called Chafee-Infante problem,
\[
u_t-\frac{1}{\lambda}u_{xx}+u^3-u=0,x\in[0,\pi] \] with boundary
conditions $u(0,t)=0,u(\pi,t)=0$.
We used 201 finite difference discretization points.
\begin{figure}[htbp]
\begin{center}
\includegraphics[width=5in]{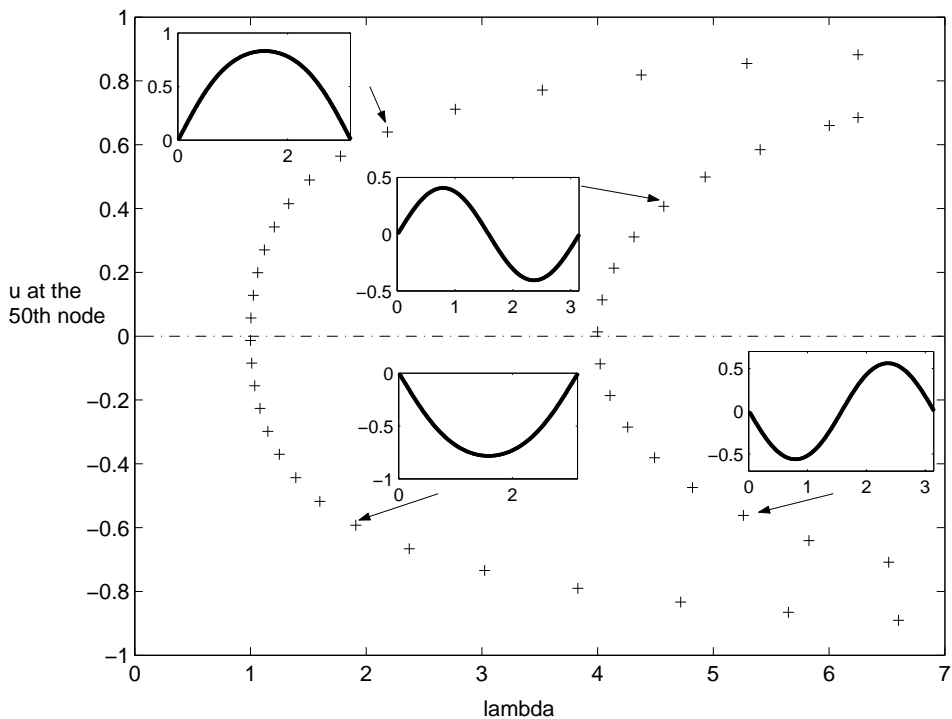}
\caption{Bifurcation diagram for the Chafee-Infante reaction
diffusion problem.} \label{kbdiag}
\end{center}
\end{figure}

\begin{figure}[htbp]
\begin{center}
\includegraphics[width=4.5in]{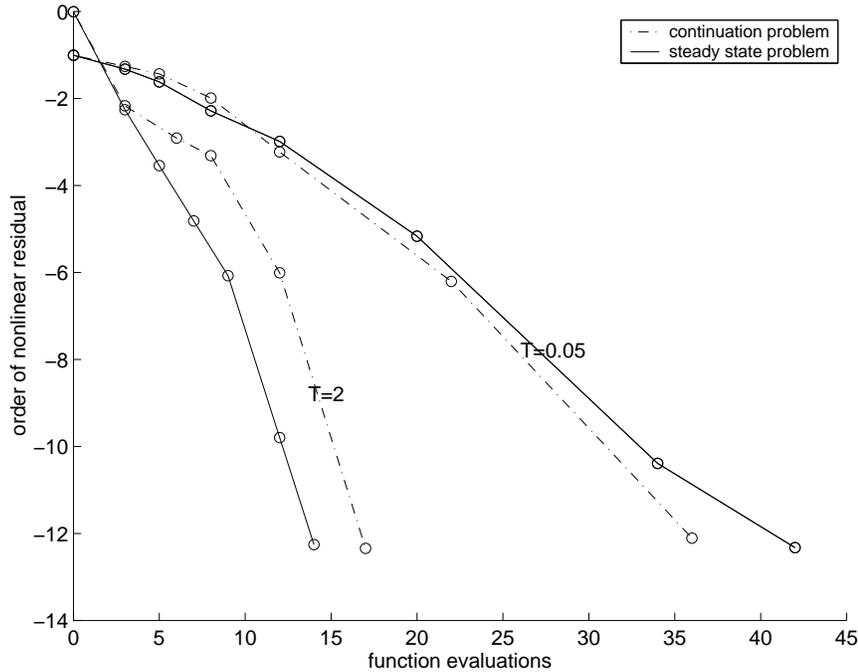}
\end{center}
\caption{Convergence of time-stepper based fixed point computation for the steady state and
the continuation of the discretized RD problem. Every circle corresponds to a Newton step. A
single ``function evaluation" consists of integration over one time
reporting horizon. } \label{twosteppers}
\end{figure}

Figure.~\ref{kbdiag} shows the bifurcation diagram for this
problem for a range of parameter values ($0 < \lambda <7$) where
up to five different spatially structured steady states exist.
Our computational tests were performed on the upper stable
solution branch close to $\lambda = 2.$
In this example we studied both the steady state problem (setting
the right-hand-side of the finite difference equations equal to
zero) and the time-stepper formulation (our integration routine
was ode15s).
The steady state (equations) / fixed point (time-stepper)
and continuation problems were solved with
Newton-GMRES solver \it nsoli \rm for various time reporting
horizons.
The results for  $\lambda=2.1386697$ are shown in
Fig.~\ref{twosteppers}.
Relative and absolute tolerances were chosen to be
$10^{-12}$.
The initial guess for the direct solution was chosen to be the
true solution plus a perturbation function $0.1\sin(x)$ .

One thing should be made clear at this point.
Using several time steps of an implicit integrator, with
the concomitant nonlinear solves, is clearly not an efficient way
of solving a fixed point problem (a single nonlinear solve).
The integrator is used here as a ``legacy code", a code that
one cannot modify.
It is in the context of such legacy codes that our approach
becomes useful, as well as in the case of multiscale computations,
where the time evolution is performed by a simulator at
a different level of description (e.g. Lattice-Boltzmann or
kinetic Monte Carlo).

\subsection{Clustering of Eigenvalues}


\begin{figure}[htbp]
\begin{center}
\includegraphics[width=4.5in]{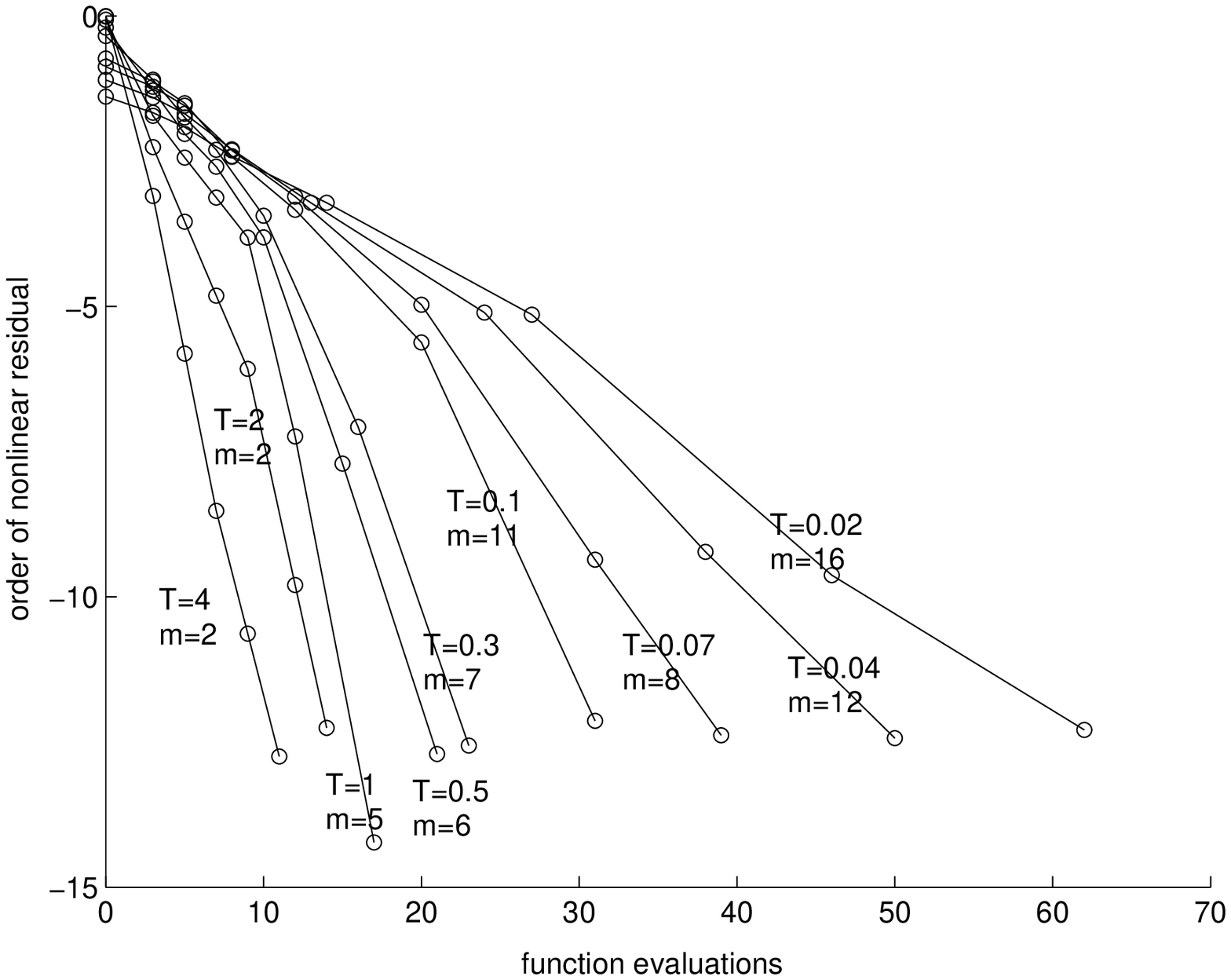}
\end{center}
\caption{Convergence of time-stepper based fixed point Newton-GMRES computation
for different time reporting horizons. Every circle corresponds to a Newton step, and
m is the number of function evaluations for the last Newton step.} \label{fig5}
\end{figure}


Using the Finite Difference Method, the ODE system from discretizing the
original PDE has a form of $\dot u(t)=f(u),u \in R^n $. At a
steady state $u^* $ , the matrix $ f_u $ has $n$ eigenvalues: $
\sigma _i,i=1,\ldots,n$. If we solve for $ u^* $ with a time-stepper, the system to be solved is $ u(0)- \Phi_T(u;\lambda)=0$,
which also has $n$ eigenvalues: $ 1-e^{\sigma _i T},i=1,\ldots,n$.

We have used the forward in time time-stepper to compute both stable
and unstable steady states. We have marked the unstable steady
state we computed using a reporting horizon of $T=0.1$ for
$\lambda=4.5710239$ using forward in time time-stepping. When the
parameter $\lambda$ is equal to 2.183867, $ u^* $ is a stable
steady state and all the eigenvalues are negative. When we
increase $T$, all the eigenvalues $ 1-e^{\sigma_i T}$ are
approaching (``clustering at") 1. Clustering of
eigenvalues is known to be beneficial for GMRES performance. Conversely, when we decrease $T$,
eigenvalues start leaving the cluster. This results in additional GMRES iterations.

To quantify the dependence of the performance of the iteration
on $T$ we use the number of GMRES iterations at the final step in
Figure~\ref{fig5}.
When $T$ is large, we consistently see $2$ linear iterations.
As
$T$ is reduced, we see an increase in GMRES iterations, as expected.
In this particular example, when $T = 1.78$, the number of GMRES
iterations increases to $3$.

\begin{figure}[htbp]
\begin{center}
\includegraphics[width=4.5in]{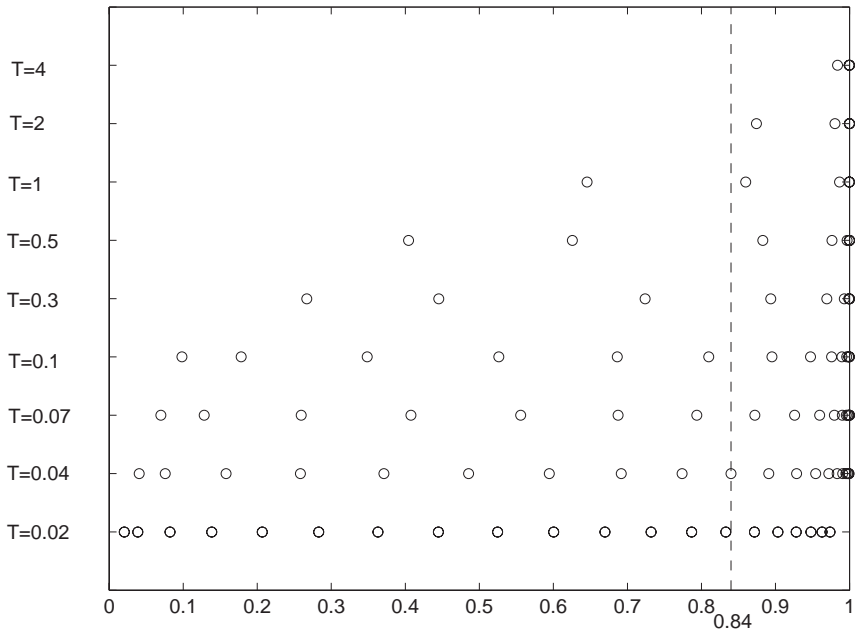}
\end{center}
\caption{Twenty smallest eigenvalues of the linearized time-stepper at
the steady state for different time reporting horizons; the dashed line corresponds to the eigenvalue
that first leaves the cluster when $T=1.78$ .}
\label{fig6}
\end{figure}

In at attempt to quantify this further we identify a ``cluster''
of eigenvalues that seems to be correlated with the performance
of the GMRES iteration.
In Fig.~\ref{fig6}, we treat those eigenvalues in the interval
[0.84 1] as ``in the cluster". The number of eigenvalues outside
the cluster (or smaller than 0.84) and the number of function
evaluations needed to finish the last Newton step are compared in
Tab.~\ref{tab1} below. A clear, strong correlation emerges.

\begin{table}[htbp]
\begin{center}
\begin{tabular}{|c|ccccccccc|} \hline
&&&&&&&&&\\
T & 4 & 2 & 1 & 0.5 & 0.3 & 0.1 & 0.07 & 0.04 & 0.02 \\
&&&&&&&&&\\ \hline

Number of &&&&&&&&&\\
eigenvalues outside & 0 & 0 & 1 & 2 & 3 & 6 &7&10&14 \\
the cluster &&&&&&&&& \\ \hline

Function evaluations &&&&&&&&& \\
needed for clustering & 2 & 2 & 2 & 2 & 2 & 2 & 2 & 2 & 2 \\
eigenvalues  &&&&&&&&& \\ \hline


Actual &&&&&&&&& \\
total function &2&2&5&6&7&11&8&12&16 \\
evaluations &&&&&&&&& \\ \hline

\end{tabular}
\end{center}
\caption{}
\label{tab1}
\end{table}

\subsection{Conclusion}
\label{sec:conclusion}  We have extended and sharpened results for
the performance of GMRES for discretizations of compact fixed
point problems to the special class of problems that arise in a
time-stepper context.
The key feature leading to the enhancement of GMRES performance is
the compactification of the spectrum. In previous work, a Green's
function based reformulation of the steady state problem for
elliptic PDEs \cite{ctk:ferng1}
led to a linearization that was a compact
perturbation of the identity, and an efficient solution via
Newton-GMRES.
Time-steppers compactify the spectrum of the linearization in a
natural way, and their properties can be exploited to obtain
accurate bounds on the convergence rates of the linear iterations
in a Newton-GMRES continuation.
We show that the additional equation in a pseudo-arclength
formulation of a parameter-dependent family of nonlinear equations
adds at most two GMRES iterations when the eigenvalues are well
separated, and obtain a bound on the convergence in terms of the
separation of the spectrum and the dimension of the slow subspace (associated with a slow/inertial
manifold for the dynamics).

We reported on numerical results that support the theory, first with
an integral equation, for which we can numerically demonstrate
that all but two of the eigenvalues of the linearization lie in a
tight cluster about $1$. The second example was a parametric study
of the steady state of a discretized parabolic partial differential equation implemented through
time-steppers.

The ``natural" compactification of the spectrum using
time-steppers provides a natural connection with the performance
of matrix-free iterative methods.
This compactification may prove useful in writing computational
wrappers, that will accelerate the convergence of legacy dynamic
simulators to stationary states.
We also expect this compactification to assist in writing
computational wrappers that will assist dynamic simulators at a
different level of model description (e.g. kinetic Monte Carlo,
Brownian Dynamics or Molecular Dynamics codes\cite{makeev,siettos,hummer}) to locate stationary states and perform continuation/bifurcation analysis of
macroscopic system observables in the so-called equation-free
framework \cite{coarse2}.

\clearpage

\bibliography{final}
\bibliographystyle{siam}

\end{document}